\documentclass{amsart}
\usepackage{amsmath,amssymb,amscd,graphicx,epsfig}
\everymath{\displaystyle}

\newtheorem{teo}{Theorem}[section]
\newtheorem{lem}[teo]{Lemma}
\newtheorem{cor}[teo]{Corollary}

\newtheorem{prop}[teo]{Proposition}
\newtheorem{lem-defi}[teo]{Lemma-Definition}

\newcommand{\mr}{\mathbb{R}}
\newcommand{\mt}{\mathbb{T}}

\newcommand{\mg}{\mathbb{G}}
\newcommand{\mz}{\mathbb{Z}}

\newcommand{\mn}{\mathbb{N}}
\newcommand{\mo}{\mathbb{O}}
\newcommand{\mq}{\mathbb{Q}}

\newtheorem{theo}{Theorem}[section]

\newtheorem{lemma}[theo]{Lemma}

\newcommand{\Aa}{{\mathcal A}}

\newcommand{\Cc}{{\mathcal C}}
\newcommand{\Dd}{{\mathcal D}}

\newcommand{\Ff}{{\mathcal F}}

\newcommand{\Ll}{{\mathcal L}}

\newcommand{\Pp}{{\mathcal P}}
\newcommand{\Rr}{{\mathcal R}}

\newcommand{\Vv}{{\mathcal V}}

\newcommand{\ZM}{{\mathbb Z}}

\begin{document}
\title {rotation topological factors of minimal $\ZM^{d}$-actions on
the cantor set}

\pagestyle{headings}
\author{Maria Isabel Cortez, Jean-Marc Gambaudo and Alejandro Maass}
\address{{\it M. I. Cortez }: Departamento de Ingenier\'{\i}a
Matem\'atica, Fac. Ciencias F\'{\i}sicas y Ma\-te\-m\'a\-ti\-cas,
Universidad de Chile,
Av. Blanco Encalada 2120 5to piso, Santiago, Chile \\ and Institut
de Math\'ematiques de Bourgogne, U.M.R. CNRS 5584, Universit\'e de
Bourgogne, U.F.R. des Sciences et T\'echniques, B.P. 47870-\,\,
21078 Dijon Cedex, France }
\email{mcortez@dim.uchile.cl}

\address{{\it J.-M. Gambaudo}:Centro de Modelamiento Matem\'atico,
U.M.R. CNRS 2071,
Av. Blanco Encalada 2120,  7to piso, Santiago, Chile.}
\email{gambaudo@u-bourgogne.fr}

\address{{\it A. Maass}:  Departamento de Ingenier\'{\i}a Matem\'atica and
Centro de Modelamiento Matem\'atico,
Fac. Ciencias F\'{\i}sicas y Matem\'aticas, Universidad de Chile,
Av. Blanco Encalada 2120 5to piso, Santiago, Chile}
\email{amaass@dim.uchile.cl}

\date{\today}

\begin{abstract}
In this paper we study conditions under which a free
minimal $\mz^d$-action on the Cantor set is a topological extension
of the action of $d$  rotations, either on the product $\mt^d$
of $d$ $1$-tori or on a single $1$-torus $\mt^1$. We extend
the notion of {\it linearly recurrent} systems  defined for
$\mz$-actions on the Cantor set to $\mz^d$-actions
and we derive in this more general setting, a necessary and
sufficient condition, which involves a natural combinatorial data
associated with the action, allowing the existence of a rotation
topological factor of one these two types.
\end{abstract}

\maketitle
\markboth{M.I. Cortez, J.-M. Gambaudo and A.  Maass}
{Rotation factors of minimal $\mz^d$-actions on the Cantor set}

\section{Introduction}

Let $(X,{\Aa})$ be a $\mz^d$-action (by homeomorphisms) on a
compact metric space $X$. The action is
{\it free} if ${\Aa}(\bar n,x)=x$ for some $\bar n \in
\mz^d$ and $x \in X$ implies $\bar n=0$ and  is {\it minimal} if
the orbit of any point $x \in X$, $O_{\Aa}(x)=\{\Aa(\bar n,x):
\bar n \in \mz^d \}$, is dense in $X$.

The simplest non trivial examples of free minimal
$\mz^d$-actions on a compact metric space are given by
``rotation-type'' actions on compact topological groups.
This type of factors
play a central role in topological dynamics of $\mz^d$-actions
since in particular they determine weak mixing property through
the existence of continuous eigenvalues. In this paper,
we focus on two kinds of ``rotation-type'' factors that we describe now.

\begin{itemize}
\item
First consider the
$\mz^d$-action generated by $d$ rotations on the product
$d$-torus $\mt^d =\mr^d/\mz^d =\mt^1\times \dots \times \mt^1$,
each rotation acting on $\mt^1$. More precisely, take
$\bar\theta=(\theta_1, \dots, \theta_d) \in \mr^d$ and let
${\Aa}^d_{\bar\theta}: \mz^d\times \mt^d \to \mt^d$ be the map
defined by:
$${\Aa}^d_{\bar\theta}(\bar n,x)=x + [\bar n,\bar\theta]\mod \mz^d \ ,$$
for $\bar n=(n_1, \dots , n_d) \in \mz^d$, $x \in \mt^d$ and
where $[\bar n,\, \bar\theta]\,=\,(n_1 \cdot \theta_1, \dots , n_d
\cdot \theta_d)$. This construction yields a minimal  $\mz^d$-action
$(\mo^d,  \Aa^d_{\bar\theta})$ on the closure $\mo^d$ of the orbit of $0$
in the $d$-torus $\mt^d$. When the coordinates of $\bar \theta$  are rationally independent, 
the set $\mo^d$ is the $d$-torus $\mt^d$ and the action is free.

\item The same
$\bar\theta$ can be used to define a $\mz^d$-action
on $\mt^1$.  Consider the map  $\Aa^1_{\bar\theta}: \mz^d\times \mt^1 \to
\mt^1$  given by
$$\Aa^1_{\bar\theta}(\bar n,t)=t + <\bar n,\bar\theta> \mod \mz,$$
where $<\cdot,\cdot>$ is the usual inner product in $\mr^d$. The $\mz^d$-action
$(\mo^1,  \Aa^1_{\bar\theta})$ on the closure $\mo^1$ of the orbit of $0$
in the $1$-torus $\mt^1$ is again minimal. When
the coordinates of $\bar\theta$ are independent on $\mq$,
the set $\mo^1$ is the $1$-torus $\mt^1$ and the action is free.
\end{itemize}

Assume $X$ is a Cantor set, ${\it {i.e.}}$, it has a countable basis of
closed and open (clopen) sets and has no isolated points (or
equivalently, it is a totally disconnected compact metric space
with no isolated points).
\medskip

{\it The main question we address in this paper is to determine whether
a free minimal
$\mz^d$-action $\Aa$ on the Cantor set $X$ is an extension of an
action of type $(\mo^d,  \Aa^d_{\bar\theta})$ or
$(\mo^1,  \Aa^1_{\bar\theta})$ for some $\bar\theta \in \mr^d$.}
\medskip

Notice that a complete combinatorial answer to this question is given in
\cite{BDM} in the particular case when the dimension $d =1$ and
when the free minimal $\mz$-action is {\it linearly recurrent}.
The linear recurrence of a given $\mz$-action is a property that
involves the combinatorics of {\it return times} associated with a
nested sequence of clopen   sets
(for further references on linearly recurrent $\mz$-actions see
\cite{CDHM},\cite{Du1} and \cite{Du2}).

The notion of return time to a clopen set can be generalized to
$\mz^d$-actions when $d\geq 2$. In this case, the combinatorics of the
return times associated with a nested sequence of clopen sets inherits
a richer structure than in the case $d =1$ . However, as for  $d=1$,
there exists a natural definition of linearly recurrent $\mz^d$-action.
These generalizations are developed in Section \ref{comb} which is devoted
to the combinatorics of return times (for further references on the structure of return times associated with a $\mz^d$-action see \cite{BG} where the hierarchical ideas used in this paper are introduced, see also \cite{S} and \cite{SW} for related topics).  

This combinatorial approach allows us to derive a necessary condition on
the action to be an extension of an action of one of the two rotations
described above. In the case of a linearly recurrent action this condition
is sufficient. This result is given in Section \ref{main}
(Theorem \ref{theoreme}) together with its proof.

\section{Combinatorics of  return times}
\label{comb}

Let us start this section with some general considerations.

Let $\mr^d$  be the  Euclidean  $d$-space  and $\Vert-\Vert$ its Euclidean norm.
Consider two positive numbers $r$ and $R$. An {\it $(r,R)$-Delone set}
is a  subset $\Dd$ of the  $d$-space $\mr^d$ equipped with the Euclidean norm  $\Vert-\Vert$,
which satisfies the following two properties:

\begin{itemize}
\item [(i)] {\it Uniformly Discrete:} each open ball with radius $r$
in $\mr^d$ contains at most one point in $\Dd$;
\item[(ii)] {\it Relatively Dense:} each open ball with radius $R$
contains at least one point in $\Dd$.
\end{itemize}

When  the constants $r$ and $R$ are not explicitly used, we will say
in short {\it Delone set} for an $(r,R)$-Delone set. We refer to
\cite{LP} for a more detailed approach of the theory of Delone sets.

A {\it patch } of a Delone set $\Dd$ is a finite subset of $\Dd$.
A Delone set is of {\it finite type} if for each $M>0$, there exist
only finitely many patches in $\Dd$ of diameter smaller than $M$ up to
translation.
Finally, a Delone set of finite type is {\it repetitive} if for each patch
$P$ in $ \Dd,$ there exists $M>0$ such that each ball with radius $M$ in
$\mr^d$ contains a translated copy of $P$ in $\Dd$.

Let $x$ be a point of a Delone set $\Dd$. The {\it Voronoi cell}
$\Vv_x$ associated with $x$ is the convex closed set in $\mr^d$ defined by:

$$\Vv_x = \{ y \in \mr^d \ : \  \forall \ x' \in \Dd, \
\Vert y - x \Vert \leq \Vert y - x' \Vert \ \} \ .$$

The union $\cup_{x\in \Dd} \Vv_x$ is a cover of $\mr^d$.
We say that two points $x$ and $x'$ in $\Dd$ are {\it neighbors} if
$\Vv_x\cap \Vv_{x'} \neq \emptyset$.

The set of return vectors associated with $\Dd$ is defined by:
$$\vec \Dd  = \{ x - y  \ : \  (x,y) \in \Dd \times \Dd \} \ .$$

\begin{lem}
\label{voisin}
Let $\Dd$ be a Delone set of finite type. Then, there exists a finite
collection $\vec \Ff$ of vectors in $\vec \Dd$ such that:
\begin{itemize}
\item$ \vec \Ff = -\vec \Ff$;
\item any vector in  $\vec \Dd$ is a linear combination with non negative
integer coefficients of vectors in $\vec\Ff$.
\end{itemize}
\end{lem}

\begin{proof}
When $\Dd$ is a Delone set of finite type, the set of vectors
$$\vec\Ff = \bigcup_{(x, x') \in \Dd \times \Dd,  (x,x') \ \rm{neighbors}}\ (x-x')$$
is finite, satisfies $\vec \Ff = -\vec \Ff$ and clearly any vector in
$\vec \Dd$ is a linear combination with non negative integer coefficients
of vectors in $\vec\Ff$.
\end{proof}

Given such a set $\vec \Ff$, we can define the {\it $\vec \Ff$-distance}
$d_{\vec \Ff}(x, x')$ as the minimal number of vectors in $\vec\Ff$
(counted with multiplicity) needed to write $x-x'$ for $x, x' \in \Dd$.
The {\it $\vec\Ff$-diameter} of a patch  $P$, denoted by
$diam_{\vec \Ff}(P)$, is the maximal $\vec\Ff$-distance of pair
of points in $\Dd$.
\bigskip

Consider now a free minimal $\mz^d$-action $\Aa$ on the Cantor set $X$.
Let $\Cc$ be a clopen set in $X$ and $y$ a point in $\Cc$.
The set of {\it return times} of the orbit of $y$ in $\Cc$ is defined by
$$\Rr_\Cc(y) = \{\bar n\in \mz^d \  : \ \Aa(\bar n, y) \in \Cc \} \ .$$

\begin{prop}
\label{delone}
The set of return times $\Rr_\Cc(y)$ is a repetitive Delone set of
finite type in
$\mz^d$. Furthermore, if $y$ and $y'$ are two points in $\Cc$,
the  sets $\Rr_\Cc(y)$ and $\Rr_\Cc(y')$ have the same patches up to
translation.
\end{prop}

\begin{proof}
$\bullet$ {\it  $\Rr_\Cc(y)$  is a Delone set of finite type.}

The minimality of the action implies that the orbit of any point
in $X$  visits $\Cc$. For each $x \in X$ consider
$\bar n_x \in \mz^d$ be such that $\Aa(\bar n_x, x)$
is in $\Cc$.
Since $\Cc$ is open, there exists a small neighborhood $U_x$ of
$x$ such that for any $x'$ in $U_x$ we also have $\Aa( \bar n_x, x')\in \Cc$.
Therefore $\{ U_x \ : \ x \in X \}$ is a cover of $X$.
Since $X$ is compact, we can extract a finite cover
$\{U_{x_i} \ : \ i \in I \}$.
Let us choose $R> \max_{i\in I} \Vert \bar n_{x_i}\Vert$. It is clear
that any ball with radius $R$ in $\mr^d$
intersects $\Rr_\Cc(y)$. Thus, $\Rr_\Cc(y)$ is relatively dense.
Since it is a subset of $\mz^d$, it is a Delone set of finite type.

$\bullet$ {\it  $\Rr_\Cc(y)$  is repetitive\footnote{
The proof that minimality implies repetitivity  is classical
and works  in a more general situation. However, for sake of completeness,
we fix it here for our specific context.}.}

Consider a patch $P$ in $\Rr_\Cc(y)$, choose $\bar n_0$ in $P$ and
let $z =\Aa( \bar n_0, y) \in \Cc$. Choose now a clopen set $\Cc_z$
containing $z$, small enough so that for any $z'$ in $\Cc_z$,
$\Aa(\bar n- \bar n_0, z')$ is in $\Cc$ for each $\bar n$ in $P$.
The set $\Rr_{\Cc_z}(z)$ is relatively dense, let $R_1$  be its
$R$-constant. Let $M$ stand for the diameter of $P$ and let us prove that 
any ball with radius $R_1+M$ in $\mr^d$ contains a
translation of the patch $P$.
Indeed, given such a ball $B$,
choose an element $\bar m \in \Rr_{\Cc_z}(z)$ in the corresponding
centered sub-ball of radius $R_1$, then by construction $\bar m + P$ belongs
to $\Rr_\Cc(y)$ and to the ball $B$.

$\bullet$ {\it $\Rr_\Cc(y)$ and $\Rr_\Cc(y')$ have the same patches up to
translation.}

Let $P$ be a patch of $\Rr_\Cc(y)$ and $\bar n_0$ be a point in $P$.
The minimality of the action implies that the orbit of $y'$ accumulates
on $z= \Aa( \bar n_0, y)$. This means that there exists $\bar n_1 \in \mz^d$
such that $\Aa(\bar n_1 +\bar n -\bar n_0, y')$ is in $\Cc$ when
$\bar n$ is in $P$. Thus a translation of the patch $P$ is in $\Rr_\Cc(y')$.
\end{proof}
\bigskip

The set of {\it return vectors} associated with $ \Cc$ is defined by:

$$\vec\Rr_\Cc = \Rr_\Cc(y) - \Rr_\Cc(y)  = \{\bar n - \bar m \ : \
(\bar n, \bar m)\in \Rr_\Cc(y)\times \Rr_\Cc(y)\} \ .$$

The fact that for any pair of points $y$ and $y'$ in $\Cc$,
the patches of $\Rr_\Cc(y)$ and $\Rr_\Cc(y')$ fit up to translation,
implies that $\vec\Rr_\Cc$ does not depend on $y$ in $\Cc$,
as suggested by the notation.
Lemma \ref{voisin} and Proposition \ref{delone} yield the following corollary.

\begin{cor}
\label{module}
There exists in $\vec\Rr_\Cc$ a finite collection of vectors $\vec\Ff_\Cc$
such that:
\begin{itemize}
\item $\vec\Ff_\Cc = - \vec\Ff_\Cc$;
\item any vector in $\vec\Rr_\Cc$ is a linear combination with non negative
integer coefficients of vectors in $\vec\Ff_\Cc$.
\end{itemize}
\end{cor}

Such a  set $\vec \Ff_\Cc$ is called a {\it set of first return vectors}
associated with $\Cc$.
\bigskip

Now we shall construct a combinatorial data associated to
a $\mz^d$-action. Let $x$ be a point in $X$ and consider
a sequence of nested clopen sets
$X= \Cc_0\supseteq \Cc_1\supseteq\dots\supseteq \Cc_n$ such that
$$\bigcap_{n\geq 0}\Cc_n\, =\, \{x\} \ .$$

Consider also the associated sets of return times $\Rr_{\Cc_n}(x)$,
of return vectors $\vec \Rr_{\Cc_n}$ and of first return vectors
$\vec \Ff_{\Cc_n}$ that we denote respectively (in short) by $\Rr_n(x)$,
$\vec \Rr_{n}$ and $\vec \Ff_{n}$.

\begin{prop}\label{chile}
For each $n\geq 0$, there exist a constant $k(n)>0$
and a partition of  $\Rr_n(x)$ in disjoint patches
$\{\Pp_n(\bar m)\}_{\bar m\in \Rr_{n+1}(x)}$ such that,
for each $\bar m\in \Rr_{n+1}(x)$:
\begin{itemize}
\item [(i)]  $\Pp_n(\bar m)\cap\Rr_{n+1}(x) = \{\bar m\}$;
\item  [(ii)] $diam_{\vec\Ff_n}(\Pp_n(\bar m))\leq k(n)$.
\end{itemize}
\end{prop}
\begin{proof}
For any point $\bar m$ in $\Rr_{n+1}(x)$ consider its Voronoi cell
$\Vv_{\bar m,n+1}$. The intersection  of this Voronoi cell  with $\Rr_n(x)$ defines a patch
$\Pp_n(\bar m)$ which intersects $\Rr_{n+1}(x)$ at $\bar m$.  It may occasionally happen that a point  $\bar l$ in 
$\Rr_n(x)$  belongs to more than one Voronoi cell
$\Vv_{\bar m,n+1}$. It this case, we make an arbitrary choice to exclude the point $\bar l$ from all the patches it belongs to but one.  This surgery done, the collection  of patches $\{\Pp_n(\bar m)\}_{\bar m\in \Rr_{n+1}(x)}$  realizes a partition of $\Rr_n(x)$.  Furthermore, since $\Rr_{n+1}(x)$ and $\Rr_{n}(x)$ are
repetitive Delone sets,
the Euclidean diameters of the cells $\Vv_{\bar m,n+1}$
are bounded independently of $\bar m$, and thus their
$\vec\Ff_n$-diameters are bounded independently of $\bar m$.
\end{proof}
\bigskip

The data  $(x, \{\Cc_n\}_{n\geq 0},
\{\{\Pp_n(\bar m)\}_{\bar m\in \Rr_{n+1}(x)}\}_{n\geq 0},
\{\vec\Ff_n\}_{n\geq 0})$ is called {\it a combinatorial data }
associated with the action $(X, \Aa)$.

We remark that Proposition \ref{chile} does not require any condition
on the nested sequence of clopen sets. By forgetting some $\Cc_n$'s in
the sequence, it is always possible to insure the following two extra
properties for the combinatorial data:
\begin{itemize}
\item [(iii)]  for each  $n\geq 0$ and for each $\bar m$ in $\Rr_{n+1}(x)$
$$\vec \Ff_n\,  \subseteq \, \Pp_n(\bar m)\, -\,  \Pp_n(\bar m);$$
\item [(iv)]  for each  $n\geq 0$ and for each $\bar m$ in $\Rr_{n+2}(x)$,
all the patches $\Pp_n(\bar m)$ are identical up to translation.
\end{itemize}

In this case, we say that the combinatorial data
$$(x, \{\Cc_n\}_{n\geq 0}, \{\{\Pp_n(\bar m)\}_{\bar m\in
\Rr_{n+1}(x)}\} _{n\geq 0}, \{\vec\Ff_n\}_{n\geq 0})$$
is {\it well distributed}.

Let $m$ and $n$ be two integers such that $0\leq n\leq m$,
and let $\bar p$ be a point in $\Rr_m(x)$. We denote by
$\Pp_n^m(\bar p)$ the patch in $\Rr_n(x)$ defined recursively by:

$$\Pp_{m-1}^m(\bar p)\, = \, \Pp_{m-1}(\bar p) \ , $$

and

$$ \Pp_n^m(\bar p)\, =\, \cup_{\bar q\in \Pp^m_{n+1}(\bar p)}
\Pp_n(\bar q) \ .$$
We adopt the convention $\Pp_m^m(\bar p)=\{\bar p\}$.
The proof of the following result is plain.

\begin{cor}
\label{sum}
For any $n_0\geq 0$ and any $\bar p$ in $\Rr_{n_0}(x)$, there exists a
unique $m_0 \geq n_0$ and a unique sequence
$\{\bar p_{l}\}_{0\leq l\leq m_0- n_0}$ of points in $\mz^d$ such that:
\begin{itemize}
\item $m_0$ is the smallest $m\geq n_0$ for which
$\bar p \in \Pp^m_{n_0}(0)$;
\item $\bar p_{0} = 0$;
\item $\bar p_l\in \Pp_{m_0-l}(\bar p_{l-1})$ and
$\bar p \in \Pp^{m_0-l}_{n_0}({\bar p_l})$ for all $1\leq l\leq m_0-n_0$;
\item $\bar p_{m_0- n_0} = \bar p$.
\end{itemize}
\end{cor}

When the constant $k(n)$ in Proposition \ref{chile} is bounded
independently on $n$, we say that the free minimal $\mz^d$-action
$\Aa$ on the Cantor set $X$ is {\it linearly recurrent}.
In this case, the combinatorial data
$(x, \{\Cc_n\}_{n\geq 0},
\{\{\Pp_n(\bar m)\}_{\bar m\in \Rr_{n+1}(x)}\}_{n\geq 0},
\{\vec\Ff_n\}_{n\geq 0})$ is said {\it adapted} to the action.

\section{Main results}
\label{main}

To each vector $\bar\theta$ in $\mr^d$ we associate the linear maps
$c^{1}_{\bar\theta} \in \Ll(\mz^d, \mt^1)$ and
$c^{d}_{\bar\theta} \in\Ll(\mz^d, \mt^d)$ defined by
$$c^{1}_{\bar\theta}(\bar p) \,=\, <\bar\theta,\bar p> \, mod\,
\mz \,\, \text{ and } \,\, c^{d}_{\bar\theta}(\bar p) \,=\,
[\bar\theta,\bar p] \, mod\, \mz^d$$ for each $\bar
p$ in $\mz^d$.

Consider a minimal free $\mz^d$-action $(X,\Aa)$ on the Cantor
set $X$ and a combinatorial data
$(x,\{\Cc_n\}_{n\geq 0}, \{\{\Pp_n(\bar m)\}_{\bar m\in
\Rr_{n+1}(x)}\} _{n\geq 0}, \{\vec\Ff_n\}_{n\geq 0})$
associated with this action. For any $n\geq 0$
and any $\bar\theta \in \mr^d$ we define the
{\it $\bar\theta$-length of $\vec\Ff_n$ of dimension $1$ and $d$}
respectively by:
$$l^{1}_{n,\bar\theta} = \max_{r_n\in \vec\Ff_n}\vert\vert\vert
c^{1}_{
\bar\theta}(r_n)\vert\vert\vert \,\,\, {\rm{and}}\,\,\,l^{d}_{n,\bar\theta}
\,=\, \max_{r_n\in \vec\Ff_n}\vert\vert\vert c^{d}_{
\bar\theta}(r_n)\vert\vert\vert, $$
where $\vert\vert\vert \cdot \vert\vert\vert$ stands for
the Euclidean distance to $0$ on the $k$-torus, $\mt^k = \mr^k/\mz^k$, $k = 1, \, d.$.
The following theorem is the main result of this paper.

\begin{teo}\label{theoreme}
Let $(X,\Aa)$ be a free minimal $\mz^d$-action on the Cantor
set $X$,
$(x, \{\Cc_n\}_{n\geq 0}, \{\{\Pp_n(\bar m)\}_{\bar m\in
\Rr_{n+1}(x)}\} _{n\geq 0}, \{\vec\Ff_n\}_{n\geq 0})$
be an associated  combinatorial data and  $k=1$ or $k =d$.
\begin{itemize}
\item[(i)]  Assume that for some $\bar\theta \in \mr^d$,
$(X,\Aa)$ is an extension of the action $(\mo^k,  \Aa^k_{\bar\theta})$.
Assume furthermore that the
combinatorial data
is well distributed. Then the series $\sum_{n\geq 0}
l^{k}_{n,\bar\theta}$ converges.
\item [(ii)] Conversely  assume that the action is linearly recurrent, that the combinatorial data is adapted to the action and that, for some
$\bar\theta \in \mr^d$, the series $\sum_{n\geq 0}
l^{k}_{n,\bar\theta}$ converges. Then
$(X,\Aa)$ is an extension of the action $(\mo^k,  \Aa^k_{\bar\theta})$.

\end{itemize}
\end{teo}

{\bf Remark 1:}
In the particular case when the $\mz^d$-action $\Aa$ is the
product of $d$  linearly recurrent $\mz$-actions on $X$,
Theorem \ref{theoreme} for $k =d$ is a direct corollary of its $d=1$
version proved in \cite{BDM}.

{\bf Remark 2:} The lie group structure of $\mt^k$ allows us to
construct a continuous surjective map $\phi: \mt^d\to \mt^1$
defined by $\phi(\alpha_1, \dots , \alpha_d)  = \alpha_1 + \dots +
\alpha_d.$ Assume that $h: (X,\Aa) \to (\mo^d,
\Aa^d_{\bar\theta})$ is an extension, then the map $\phi\circ h:
(X,\Aa) \to (\mo^1,  \Aa^1_{\bar\theta})$ is also an extension.
This is coherent with the fact  that the convergence of the series
$\sum_{n\geq 0} l^{d}_{n,\bar\theta}$  implies the convergence of
the series $\sum_{n\geq 0} l^{1}_{n,\bar\theta}$.

{\it Proof of Theorem \ref{theoreme}.} The proofs of both
assertions of Theorem \ref{theoreme} for $k=1$ or $k=d$ follow the
same scheme and will be gathered in a single demonstration. Let
$<<\cdot,\cdot >>$ stand for  $[\cdot,\cdot] \mod \mz^d$ when  $k=d$ and for
$<\cdot,\cdot> \mod \mz$ when  $k=1$.

$(i)$ Assume that the  free minimal $\mz^d$-action $(X,\Aa)$ is an
extension of the action $\Aa^k_{\bar\theta}$ on the closure
$\mo^k$ of the  orbit of the point $0$ in the $k$-torus $\mt^k$
for some $\bar\theta$ in $\mr^d$. Let us denote by $h:X \to \mo^k$
the  extension. Choose a well distributed associated combinatorial data
$$(x, \{\Cc_n\}_{n\geq 0}, \{\{\Pp_n(\bar m)\}_{\bar m\in \Rr_{n+1}(x)}\} _{n\geq 0},
\{\vec\Ff_n\}_{n\geq 0})$$
and fix $h(x) =0 \in \mt^k$.

For each $n \geq 0$ let $v_n$ be the first return vector in
$\vec\Ff_n$
 such that:
$$l^k_{n,\bar\theta} \,=\, \max_{u_n\in \vec\Ff_n}\vert\vert\vert
c^k_{\bar\theta}(u_n)\vert\vert\vert\,=\, \vert\vert\vert
c^k_{\bar\theta}(v_n)\vert\vert\vert.$$ The following observation
is a direct consequence of the continuity of $h$.

\begin{lemma}
\label{lem}
The quantity $l^k_{n,\bar\theta}$ goes to $0$ as $n$
goes to $\infty$. Furthermore, for each $\epsilon >0$ there exists
$N>0$ such that for any pair of  points $(\bar n, \bar m)$ in
$\Rr_N(x)\times \Rr_N(x)$, we have:
$$\vert\vert\vert h(\Aa( \bar n, x))\, -\,
h(\Aa( \bar m, x))\vert\vert\vert\,\leq\,
\epsilon.$$
\end{lemma}

Let $B$  be  the open  ball on the $k$-torus, centered at $0$, with radius  ${\sqrt k }/2$.  Fix $0<\epsilon<   {\sqrt k }/2$ and  let $N$ verifying the conclusion of
Lemma \ref{lem} for this $\epsilon$ and such that $l_{n,\bar \theta}^k \leq \epsilon$ for $ N\leq n$.
 The ball
$B$ is decomposed in $2^k$ sectors $S_{\epsilon_1, \dots,
\epsilon_k }$ with $\epsilon_i \in \{-1,1\}$ for $i\in
\{1,\dots,k\}$ defined by

$$S_{\epsilon_1, \dots, \epsilon_k } = \{(x_1, \dots, x_k) \in B \ :
\ x_i \cdot \epsilon_i\geq 0 \ , \ \forall i\in \{1, \dots, k\} \}
\ .
$$

Let $I_{\epsilon_1, \dots, \epsilon_k}$ be  the set of integers
$n$ such that $c^k_{\bar\theta}(v_n)$ is in $S_{\epsilon_1, \dots,
\epsilon_k }$ and let us prove that the series $\sum_{n\in
I_{\epsilon_1, \dots, \epsilon_k }}l^k_{n,\bar\theta}$ converges.
Actually, we only need to prove that the series $\sum_{n\in I_{1,
\dots, 1 }}l^k_{n,\bar\theta}$ converges, a similar proof works
for the other cases. This sum can be spitted into two parts:

$$\sum_{n\in I_{1, \dots, 1 }}l^k_{n,\bar\theta}\, =\,
\sum_{n \in I_{1, \dots, 1 }, \,\rm{even}}l^k_{n,\bar\theta} \, +\,
\sum_{n\in I_{1, \dots, 1 },\,\rm{odd}}l^k_{n,\bar\theta}.$$

Here again we only need to prove that the series $\sum_{n\in I_{1,
\dots, 1 }, \,\rm{even}}l^k_{n,\bar\theta}$ converges, a similar
proof works also for the case where $n$ is odd. Observe that we
are assuming $I_{1, \dots, 1 }$ is infinite.

The proof splits in five steps:

{\it Step 1:} Fix an even integer  $N_0 $ big enough  in $I_{1,
\dots, 1 }$, and let $N<n_l<n_{l-1} <\dots <n_1 <N_0$ be the
ordered sequence of even integers bigger than $N$ that belong to
$I_{1, \dots, 1 }$.

{\it Step 2:}
Consider two points $\bar m_1$ and $\bar p_1$ in $\Rr_{n_1 }(x)$ such that
$$v_{n_1} = \bar p_1 - \bar m_1.$$
Since the combinatorial data is well distributed, the two patches
$\Pp_{n_2}(\bar m_1)$ and
$\Pp_{n_2}(\bar p_1)$ are identical up to translation and
there exists a pair of points
$(\bar m_2, \bar m'_2)$ in $\Pp_{n_2}(\bar m_1)\times \Pp_{n_2}(\bar m_1)$ such that
 $$v_{n_2} = \bar m'_2  - \bar m_2.$$
We define $\bar p_2 $ in $\Pp_{n_2}(\bar p_1)$ by
$\bar p_2 -\bar p_1= \bar m_2-\bar m_1 + v_{n_2}$.  We have:
$$\bar p_2 -\bar m_2\, =\, v_{n_1} + v_{n_2}.$$

{\it Step 3:} Since the combinatorial data is well distributed,
the two patches $\Pp_{n_3}(\bar m_2)$ and $\Pp_{n_3}(\bar p_2)$
are identical up to translation and there exists a pair of points
$(\bar m_3, \bar m'_3)$ in $\Pp_{n_3}(\bar m_3)\times
\Pp_{n_2}(\bar m_2)$ such that

$$v_{n_3} = \bar m'_3  - \bar m_3.$$
We define $\bar p_3 $ in $\Pp_{n_2}(\bar p_2)$ by $\bar p_3 -\bar
p_2= \bar m_3-\bar m_2+ v_{n_3}$. We have:
$$\bar p_3 -\bar m_3\, =\, v_{n_1} + v_{n_2}+ v_{n_3}.$$

{\it Step 4:} We iterate this construction until we get the points
$\bar m_l$ and $\bar p_l$ which satisfy:
$$\bar p_l -\bar m_l\, =\, \sum_{j=1}^{l} v_{n_j}.$$

{\it Step 5:}  We have:
\begin{align*}
\vert\vert\vert h(\Aa(\bar p_l, x)) \, - \, h(\Aa(\bar
m_l,x))\vert\vert\vert\ &= \vert\vert\vert <<\sum_{j=1}^l
v_{n_j},\bar \theta>>
\vert\vert\vert \\
&=\vert\vert\vert \sum_{j=1}^{l} c^k_{\bar\theta}(v_{n_j})
\vert\vert\vert
\end{align*}
Since $\bar p_l$ and $\bar m_l$ are in $\Rr_N(x)$, Lemma \ref{lem}
implies that:
$$\vert\vert\vert
\sum_{j=1}^{l} c^k_{\bar\theta}(v_{n_j}) \vert\vert\vert \leq
\epsilon \ .$$

Let $\pi: B\to B'$ be  the canonical isometric identification of the ball $B$ with the   open  ball  $B'$ in  the Euclidean space $\mr^d$ centered at $0$ with radius  ${\sqrt k}/2$. Through this identification, it is clear that for all $x$ in $B$: $\,\vert\vert\vert x\vert\vert\vert = \vert \vert\pi(x)\vert\vert $. Moreover for any pair of points $x,\, x'$ in $S_{1, \dots, 1}$ such that  $x+x'$   is also in $S_{1, \dots, 1}$, we have: 
$\pi(x+x') = \pi(x) +\pi(x')$.  It follows that 
$$\vert\vert\vert
\sum_{j=1}^{l} c^k_{\bar\theta}(v_{n_j}) \vert\vert\vert\, =\, 
\vert\vert
\sum_{j=1}^{l} \pi(c^k_{\bar\theta}(v_{n_j}) )\vert\vert.$$

Finally, since for $1\leq j \leq l$, $c_\theta^k(v_{n_j}) $ is in $S_{1, \dots, 1}$, we have:
$$
\sum_{j=1}^{l}\vert\vert \pi(c^k_{\bar\theta}(v_{n_j}) )\vert\vert\,\leq\, 1/{\sqrt k}\cdot\vert\vert
\sum_{j=1}^{l} \pi(c^k_{\bar\theta}(v_{n_j}) )\vert\vert,$$
 which implies
$$\sum_{N\leq n, \, n \in I_{1, \dots, 1 }, \, \rm{even}}
l^k_{n_j,\bar\theta} \, \leq \, 1/{\sqrt k} \cdot \epsilon \ .$$

This insures that the series $\sum_{n\in I_{1, \dots, 1 }, \, \,
\rm{even}}l^k_{n,\bar\theta}$ converges, and consequently the
series $\sum_{n\geq 0} l^k_{n,\bar\theta}$ converges too.

\bigskip

(ii)  Let $(X,\Aa)$ be a linearly recurrent $\mz^d$-action on the
Cantor set $X$.  Assume that the combinatorial data is adapted to
the action and that the series of $\bar\theta$-lengths
$\sum_{n\geq 0} l^k_{n,\bar\theta}$ converges for some
$\bar\theta$ in $\mr^d$. Fix $\epsilon >0$ and choose $n_0 \in
\mn$ big enough so that
$$ \sum_{n\geq n_0} l^k_{n,\bar\theta} < \epsilon.$$

Let us define the map $h$ on the $\mz^d$-orbit of $ x$ by,
$$h(\Aa(\bar n, x))\, = <<\bar n, \bar\theta>> = \,\Aa^k_{\bar\theta}
(\bar n,0)$$ for each $\bar n$ in $\mz^d$. In order to prove that
the map $h$ extends to a continuous map on the closure $\mo^k$of
the orbit of $0$ in $\mt^k$, it is enough to prove that $h$ is
uniformly continuous, which follows from the continuity of $h$ at
$x$. Consider a point $\bar p$ in $\Rr_{n_0}(x)$ and  apply
Corollary \ref{sum}. There exists a unique $m_0 \geq n_0$ and a
unique sequence $\{\bar p_l\}_{0\leq l\leq m_0- n_0}$ of points in
$\mz^d$ such that:
\begin{itemize}
\item $m_0 $ is the smallest $m\geq n_0$ for which $\bar p \in \Pp^m_{n_0}(0)$;
\item $\bar p_{0} = 0$;
\item $\bar p_l\in \Pp_{m_0-l}(\bar p_{l-1})$ and $ \bar p \in \Pp^{m_0-l}_{n_0}({\bar p_l}), \
\forall \,\,1\leq l\leq m_0-n_0$;
\item $\bar p_{m_0- n_0} = \bar p$.
\end{itemize}

Let us write:
$$h(\Aa(\bar p, x)) = \sum_{l=1}^{m_0-n_0}(h(\Aa( \bar p_{l}, x) )-h(\Aa( \bar p_{l-1}, x))).$$

For any  ${1\leq l\leq m_0- n_0}$ both points $\bar p_l$ and $\bar
p_{l-1}$ are in $ \Pp_{m_0-l}(\bar p_{l-1})$. Consequently there
exists  a collection $\{v_{{m_0}-l, i}\}_{1\leq i \leq
q({m_0}-l)}$ of  vectors in $\vec \Ff_{{m_0}-l}(\bar p_{l-1})$
such that:
\begin{itemize}
\item  $q({m_0}-l)\leq k(m_0 - l)$;
\item the sequence of points $\{\bar p_{l-1, i} \}_{0\leq i\leq q({m_0}-l)}$ defined by:
{\begin{itemize}
\item $\bar p_{l-1, 0} =  \bar p_{l-1}$;
\item $\bar p_{l-1, i}  = \bar p_{l-1, i-1} + v_{{m_0}-l, i}$  for $1\leq i\leq q({m_0}-l)$;
\item $\bar p_{l-1, q({m_0}-l)} = \bar p_l,$
\end{itemize} }
belongs to $\Rr_{m_0-l}(x)$.
\end{itemize}
This yields
$$h(\Aa(\bar p, x))\, =\, \sum_{l=1}^{m_0-n_0}\sum_{i=1}^{p({m_0}-l)}
(h(\Aa(\bar p_{l-1,i}, x  )) )-h(\Aa(\bar p_{l-1, i-1}, x ))) .$$

Now we use the fact that the action is linearly recurrent and that
the combinatorial data is adapted to this action.  We denote by
$L$ a uniform upper bound for the sequence $\{k(n)\}_{n\geq 0}$.
We get,

$$\vert\vert\vert h(\Aa(\bar p, x))\vert\vert\vert \, \leq \, L
\cdot \sum_{l=1}^{m_0-n_0} l^k_{m_0-l, \bar\theta}\, \leq \, L
\cdot \sum_{n=n_0}^{\infty} l^k_{n, \bar\theta}\,\leq\,
\epsilon.$$

This proves the continuity of $h$ at $x$.
\hfill$\Box$
\bigskip

{\bf Acknowledgments.}  All authors acknowledge financial support
from ECOS-Conicyt grant C03-E03. The third author thanks also
support from Programa Iniciativa Cient\'{\i}fica Milenio P01-005
and FONDECYT 1010447.

\enddocument